\documentclass[11pt,a4paper,sunlimits]{article}
\usepackage{amsmath,amsfonts}
\usepackage{amsthm}
\usepackage{amssymb}
\usepackage{amscd}
\usepackage{xypic}
\usepackage{verbatim}
\usepackage[plainpages=false,colorlinks,hyperindex,pdfpagemode=None,bookmarksopen,linkcolor=red,citecolor=blue,urlcolor=blue]{hyperref}
\usepackage{pdflscape}
\usepackage{stmaryrd}
\usepackage{mathabx}
\usepackage{multirow}
\usepackage{appendix}
 \swapnumbers
\theoremstyle{plain}

\renewcommand{\marginpar}[1] {  }
\renewcommand{\comment}[1] {  }

\usepackage[all]{xy}
\catcode`\é=\active\def é{\'e}
\catcode`\è=\active\def è{\`e}
\catcode`\à=\active\def à{\`a}
\catcode`\ù=\active\def ù{\`u}

\catcode`\ê=\active\def ê{\^{e}}
\catcode`\î=\active\def î{\^{i}}
\catcode`\ô=\active\def ô{\^{o}}
\catcode`\û=\active\def û{\^{u}}
\catcode`\ç=\active\def ç{\c{c}}

\newtheorem{theo}{Theorem}[section]
\newtheorem{lem}[theo]{Lemma}
\newtheorem{prop}[theo]{Proposition}
\newtheorem{cor}[theo]{Corollary}
\newtheorem{defi}[theo]{Definition}

\theoremstyle{remark}

\numberwithin{equation}{section}

\def \g{\mathfrak}

\def\bb{\backslash}
\def\al{\alpha}
\def\be{\beta}
\def\De{\Delta}
\def\de{\delta}
\def\ep{\varepsilon}

\def\ga{\gamma}
\def\Ga{\Gamma}

\def\la{\lambda}
\def\om{\omega}
\def\Om{\Omega}

\def\si{\sigma}

\def\Te{\Theta}
\def\te{\theta}

\def\bb{\backslash}

\def\cC{{\mathcal C}}
\def\cD{{\mathcal D}}

\def\cF{{\mathcal F}}

\def\cH{{\mathcal H}}

\def\cM{{\mathcal M}}

\def\cV{{\mathcal V}}
\def\cW{{\mathcal W}}

\def\cY{{\mathcal Y}}
\def\cZ{{\mathcal Z}}
\DeclareFontFamily{OT1}{rsfs}{}
\DeclareFontShape{OT1}{rsfs}{n}{it}{<-> rsfs10}{}
\DeclareMathAlphabet{\mathscr}{OT1}{rsfs}{n}{it}

\def\Er{{\mathscr E}}

\def\A{\mathbb A}
\def\F{{\rm F}}%{\mathbb F}
\def\K{\mathbb K}
\def\C{\mathbb C}
\def\D{\mathbb D}
\def\Q{\mathbb Q}
\def\R{\mathbb R}
\def\H{\mathbb H}
\def\X{\mathbb X}
\def\Y{\mathbb Y}

\def\N{\mathbb N}
\def\Q{\mathbb Q}
\def\R{\mathbb R}
\def\S{\mathbb S}
\def\Z{\mathbb Z}

\def\bG{\mathbf G}

\def\me{\medskip}
\def\no{\noindent}
\def\dis{\displaystyle}

\def\ste{\par\smallskip\noindent}
\def\ste{\par\smallskip\noindent}

\def\dem{ {\em Proof~: \ste }}
 \def\beq{\begin{equation}}
\def\eeq{\end{equation}}

\newenvironment{res}
               {\begin{equation}\begin{minipage}{0.85\textwidth}}
               {\end{minipage}\end{equation}}
\def\ber{\begin{res}}
\def\eer{\end{res}}

\def\qed{{\null\hfill\ \raise3pt\hbox{\framebox[0.1in]{}}\break\null}}

\begin{document}

\author{P. Delorme\thanks{The first author was supported by a grant of Agence Nationale de la Recherche with reference
ANR-13-BS01-0012 FERPLAY.},   P. Harinck}
\title{Relative trace formula for compact quotient and pseudocoefficients for relative discrete series.}
\date{}
\maketitle

\begin{abstract} We introduce the notion of relative pseudocoefficient for relative discrete series of real spherical homogeneous spaces of reductive groups. We prove that such relative pseudocoefficient does not exist for semisimple symmetric spaces of type $G_\C/G_\R$ and construct strong relative pseudocoefficients for some hyperbolic spaces. We establish a toy model for the relative trace formula of H.~Jacquet for compact discrete quotient $\Ga\bb G$. This allows us to prove  that a relative discrete series which admits strong pseudocoefficient with sufficiently small support      occurs in the spectral decomposition of $L^2(\Ga\bb G)$ with a nonzero  period.   \end{abstract}
 \medskip

\section{Introduction}
It  is shown  in (\cite{CD}  Corollary of Proposition 4)      that pseudocoefficients  of all discrete series of   the group of  real points of  any  connected reductive  algebraic group defined over $\R$, always exist.

In this article we define the notion of relative pseudocoefficient for relative discrete series of real spherical homogeneous spaces of reductive groups  (see Definition \ref{pseudocoeff}).

In the case of semisimple symmetric spaces of type $G_\C/G_\R$, we show,  using properties of orbital integrals (see \cite{Bou}  and \cite{Ha3}), that  for all relative discrete series, it does not exist relative pseudocoefficients (see Theorem \ref{nopseudo}).\me

We look also to  hyperbolic spaces over $\C$ and over the quaternions $\H$, i.e $\X=G/H$, where $G=U(p,q,\K)$ and $H=U(1,\K)\times U(p-1,q,\K)$ for $\K=\C$ or $\H$.

 For the rest of this introduction we fix  such a $G$ and $H$. Using results of J. Faraut,    M.~Flensted-Jensen and K.~Okamoto  and results of \cite{DFJ}, we show  that some relative discrete  series of these spaces  do not have relative pseudocoefficients (see Proposition \ref{nopseudo}), but there exists a countable  family 
of relative discrete series having what we call strong relative pseudocoefficients i.e. which ‘’isolate ‘’ a single relative discrete series among the irreducible unitary representations of  the group $G$ having a nonzero $H$-fixed distribution vector (see Theorem \ref{Result}).

We give  an  application of the existence of these strong relative pseudocoefficients, which are of arbitrary support,  to existence of some representations occurring  in the spectral decomposition of $L^2( \Gamma\bb G)$  where $\Gamma$ is a torsion free cocompact discrete subgroup of G, stable by the involution $\sigma$ whose fixed point group is $H$. Existence of such cocompact discrete subgroups is shown  using adelic methods (see Proposition \ref{taustable}).

For this application, we establish a toy model for the more sophisticated  relative trace formula of H.~Jacquet (\cite{J}).\me

Let $\xi_\Gamma$  be the distribution vector for the right representation of $G$ in $L^2(\Ga\bb G)$ given by the integration over $\Ga\cap H\bb H$, which is compact in our case. We denote by  $c_{\xi_\Gamma, \xi_\Gamma}
$   the corresponding generalized matrix coefficient.

Then  it is easy to give   two expressions  of $c_{\xi_\Gamma, \xi_\Gamma}(f)$  for $f \in C_c^\infty (\X)$.  One,  which is called spectral, involves periods of representations, i.e. $H$-fixed distribution vectors of irreducible subrepresentations of  $L^2( \Gamma\bb G)$, the other, called geometric, involving relative orbital integrals, i.e. the average of $f$ on orbits of elements of $\Ga$ under the action of $H\times H$. We call the equality of these 2 expressions a relative trace formula.

When one plug into the spectral side  of  the relative trace formula, a strong relative pseudocoefficient, it singles out the contribution of the relative discrete series. On the geometric side, as  we can find strong relative pseudocoefficient  with sufficiently small support, we get only the contribution of the neutral element and this contribution is  equal to $f(1)$.

Altogether, it shows that the relative discrete series with strong relative pseudocoefficients    occurs in  $L^2( \Gamma\bb G)$.\me

This paper is organized as follows: In section \ref{SRTF} we prove a relative trace formula for $\Ga\bb G$ when $H$ is a unimodular closed subgroup of $G$, $\Ga$ is  a cocompat discret subgroup of $G$ such that the volume of $\Ga\cap H\bb H$ is  finite   and the centralizers of elements of $\Ga$ in $H\times H$ are unimodular. In section \ref{PSC}, we introduce the notion of relative pseudocoefficient and prove our application of existence of strong  relative pseudocoefficient when $H$ is the fixed points group of an involution of $G$. In section \ref{Strong}, we explain our results (existence or non existence of strong relative  pseudocoefficient for  relative discrete series) for hyperbolic spaces. In section \ref{exis}  we construct  cocompact discrete subgroups of unitary groups on $\C$ and $\H$ satisfying our assumptions   and the section \ref{GC} is devoted to prove that there is no relative pseudocoefficient for all relative  discrete series of a symmetric space of type $G_\C/G_\R$.\me

\no{\bf Acknowledgments}.  We thank warmly Roberto Miatello for useful discussions when elaborating  the material of section \ref{SRTF} .\\
  We thank warmly  Rapha\"el  Beuzart-Plessis and Jean-Pierre  Labesse for  their kind help  on   construction of cocompact discrete subgroups.

\section{A relative trace formula for $\Ga\bb G$.}\label{SRTF}
If $M$ is a differentiable manifold,  then $C(M)$ and  $C^\infty(M)$ will denote the space of continuous functions and  smooth  functions on $M$ respectively. Let $C_c(M)$ and $C_c^\infty(M)$ be the subspaces of compactly supported functions in  $C(M)$ and $C^\infty(M)$ respectively.\me

Let $G$ be a real reductive group. We consider a discrete cocompact subgroup $\Ga$ of $G$ and a closed subgroup $H$. We assume that
\ber\label{hyp}\begin{enumerate} \item $H$ is unimodular,
\item the volume $vol( (H\cap \Ga)\bb H)$ of $(H\cap \Ga)\bb H$  is finite,
\item for each $\ga\in\Ga$, the subgroup $(H\times H)_\ga=\{(h,h')\in H\times H; h^{-1}\ga h'=\ga\}$ is unimodular.
\end{enumerate}\eer
   We set $\Ga_H:=H\cap \Ga$. We fix Haar measure on these groups, discrete groups being equipped with the counting measure.\me
   
\no If $V$ is a topological vector space then $V'$ will denote its topological dual.\\
Let  $(\pi,V)$ a continuous representation of  $G$ in a Hilbert space $V$. We   denote by  $V^\infty\subset V$  the space of $C^\infty$ vectors  of $\pi$ endowed with its natural topology. We define the space of distribution vectors $V^{-\infty}$ as the topological dual  of  $V^\infty$. Let  $\pi_\infty$ be the representation  of $G$ in $V^\infty$, and $\pi_{-\infty}$ be the dual representation of $\pi_\infty$ in $V^{-\infty}$.\\
If $f\in C_c^\infty(G)$ and  $\xi\in V^{-\infty}$, we have  $\pi_{-\infty}(f)\xi\in (V')^\infty$. Hence, if $\xi'\in (V')^{-\infty}$,  we can define the generalized matrix coefficient $m_{\xi,\xi'}$   by
$$m_{\xi,\xi'}(f)=\langle \pi_{-\infty}(f)\xi, \xi'\rangle, \quad f\in C_c^\infty(G).$$

If $(\pi, V)$ is unitary for a scalar product $(\cdot,\cdot)$, then the map $j: v\mapsto (\cdot, v)$, intertwines the conjugate representation $(\overline{\pi},\overline{V})$ of $(\pi, V)$ and its dual representation  $(\pi',V')$. Let   $\xi\in V^{-\infty}$. We  define $\overline{\xi}\in \overline{V}^{-\infty}$ by $\overline{\xi}(w)=\overline{\xi(w)}$. By the above identification, we can consider $\overline{\xi}$ as an element of  $(V')^{-\infty}$. Thus we can  define the generalized matrix coefficient $c_{\xi,\xi}$ associated to $\xi$ by
\beq\label{coef}c_{\xi,\xi}(f)=m_{\xi,\overline{\xi}}(f)=(\pi_{-\infty}(f)\xi,\xi).\eeq \\
%%%%%%%%%%%%%%%%%%%%%%%%%%%%%%%%%%%%%
We consider the right regular representation $R$ of $G$ in $L^2(\Ga\bb G)$. \\ %%%%%%%%%%%%%%%%%%%%%%%%%%%%%%
  %%%%%%%%%%%%%%%%%%%%%%%%%%%%%%%%%%%%%%%%%%%%%%%%
%%%%%%%%%%%%%%%Then, for 
  Then, for $f\in C_c(G)$,  the corresponding operator  $R(f)$  maps any function $\psi\in L^2(\Ga\bb G)$ to the function
  $$[R(f)\psi](x)=\int_G f(g)\psi(xg) dg= \int_G f(x^{-1}y) \psi(y) dy.$$
  
  $$= \int_{\Ga\bb G}\psi(y)K_f(x,y) \; d\dot{y}$$
 where   \beq\label{Ker}\dis K_f(x,y) :=\sum_{\ga\in \Ga} f(x^{-1}\ga y),\quad x,y\in  G,\eeq\\
and this sum has only finite number of nonzero terms for $x,y$ contained in a compact subset of $G$ since $f$ is compactly supported and $\Ga$ is discrete. 

Therefore, $R(f)$ is an operator with continuous kernel $K_f$.\me
  
\no We defined the $H$-invariant linear form $\xi_\Ga$ on $C(\Ga\bb G)$,  which contains $L^2(\Ga\bb G)^\infty=C^\infty(\Ga\bb G)$, by
$$\xi_\Ga(\psi)=\int_{ \Ga_H\bb H} \psi(h) d\dot{h}.$$\\
Then,  the generalized matrix coefficient 
$$c_{\xi_\Ga,\xi_\Ga}(f):=(R_{-\infty}(f)\xi_\Ga, \xi_\Ga),\quad f\in C_c^\infty(G) $$ associated to $\xi_\Ga$ according to (\ref{coef}) is an $H$-biinvariant distribution  on $G$.\me

\no The relative trace formula in this context gives two expressions of the distribution $c_{\xi_\Ga,\xi_\Ga}$, 
the first one,   called the spectral side, in terms of irreducible representations of $G$, and the second one,  called the geometric side, in terms of orbital integrals. \me

We first deal with the spectral part. For this purpose, we consider the spectral decomposition of $L^2(\Ga\bb G)$:

$$L^2(\Ga\bb G)=\oplus_{\pi\in \hat{G}} \cH_\pi\otimes \cM_\pi,$$
where $\hat{G}$ is the set of equivalent classes of irreducible unitary representations $(\pi,\cH_\pi)$ of $G$ and $\cM_\pi$ is a finite dimensional vector space whose dimension is the multiplicity of $\pi$ in $L^2(\Ga\bb G)$.%%%%%%%%%%%%%%%%%%%%%%%%%%%%%%%%%%%
\\ Then the space $V_\pi:= \cH_\pi\otimes M_\pi$ is the $\pi$- isotypic component of $L^2(\Ga\bb G)$. We denote by $\xi_{\Ga,\pi}$ the restriction of $\xi_\Ga$ to $V_\pi$. Therefore, we obtain
\beq\label{spect}c_{\xi_\Ga,\xi_\Ga}(f) =\sum_{\pi\in \hat{G}} c_{\xi_{\Ga,\pi},\xi_{\Ga,\pi}}(f),\quad f\in C_c^\infty(G).\eeq

\no  For $\ga\in \Ga$, we   define the groups
$$(H\times H)_{\gamma}=\{(h_1,h_2)\in H\times H;\; h_1^{-1}\gamma h_2= \gamma \},\quad \textrm{and}\quad   ( \Ga_H \times  \Ga_H )_{\gamma}= (H\times H)_{\gamma} \cap (\Gamma\times \Gamma ).$$

%%%%%%%%%%%%%%%%%%%%%%%%%%%%%%%%%%%%%%%%
\begin{prop}\label{propRTF}

\begin{enumerate}  \item For   $\ga\in \Ga$, the quotient $(\Ga_H\times \Ga_H)_\ga\bb (H\times H)_\ga$ is of finite volume and for $f\in C_c^\infty(G)$, the orbital integral of $f$ at $\ga$
$$I(f, \gamma):= \int_{(H\times H)_{\gamma} \bb (H\times H)} f( h_1^{-1}\gamma h_2)dh_1\; dh_2$$
is absolutely convergent. 
\item We have the following relative trace formula
\beq\label{RTF} \sum_{\ga\in \Ga_H\bb \Ga/\Ga_H} vol((\Ga_H\times \Ga_H)_\ga\bb (H\times H)_\ga)\;I(f,\ga)=\sum_{\pi\in \hat{G}}( \pi(f)\xi_{\Ga,\pi},\xi_{\Ga,\pi}),\eeq\\
where the left hand side is absolutely convergent.
\end{enumerate}

\end{prop}
\dem The right hand side of (\ref{RTF}) is just the expression of $c_{\xi_\Ga,\xi_\Ga}(f)$ given in (\ref{spect}). \\
For the geometric side, we first express  $c_{\xi_\Ga,\xi_\Ga}(f)$ in terms of the kernel $K_f$. \me

 Let $f\in C_c(G)$. We define $\check{f}$ by $\check{f}(x)=f(x^{-1})$. 
Then, for $\psi\in C( \Ga\bb G)$, we have
$$\langle R_{-\infty}(f)\xi_\Ga, \psi\rangle=\langle\xi_\Ga, R(\check{f})\psi\rangle=\int_{\Ga_H\bb H}\big(\int_{\Ga\bb G}\psi(y) K_{\check{f}}(h,y) d\dot{y}\big)d\dot{h}.$$\\
%%%%%%%%%%%%%%%%%%%%%%%%%%
The kernel  $K_{\check{f}}$ is continuous and $|K_{\check{f}}|\leq K_{|\check{f}|}$. Applying the above equality to $|\psi|$ and 
$|f|$, one sees that the double integral on the right side is absolutely convergent and we can apply Fubini Theorem. Thus we obtain
$$\langle R_{-\infty}(f)\xi_\Ga, \psi\rangle=\int_{\Ga\bb G}\psi(y)\big(\int_{\Ga_H\bb H}K_{\check{f}}(h,y)d\dot{h}\big)d\dot{y}.$$
We deduce that     $R_{-\infty}(f)\xi_\Ga$ is the continuous function on $\Ga\bb G$ given by   $$ \big(R_{-\infty}(f)\xi_\Ga\big)(y)=\int_{\Ga_H\bb H}\overline{K_{\check{f}}(h,y)}d\dot{h}.$$\\
Therefore, we can extend the map $\varphi\in C_c^\infty(G)\mapsto c_{\xi_\Ga,\xi_\Ga}(\varphi)=(R_{-\infty}(\varphi)\xi_\Ga,\xi_\Ga)$  to $  C_c(G)$. Since $K_{\check{f}}(x,y)=K_{f}(y,x)$, we obtain  for $f\in C_c(G)$  
\beq\label{xigam}c_{\xi_\Ga,\xi_\Ga}(f)=\int_{\Ga_H\bb H}\big(\int_{\Ga_H\bb H} K_f(x,y) d\dot{x}\big)d\dot{y}= \int_{\Ga_H\bb H} \big(\int_{\Ga_H\bb H}\sum_{\ga\in\Ga} f( h_1^{-1}\gamma h_2)dh_1 \big) \; dh_2.\eeq\\
For $(h_1,h_2)\in H\times H$, we have 

\beq\label{sum}\sum_{\ga\in \Ga} f( h_1^{-1}\gamma h_2)=\sum_{[\ga]\in \Ga_H\bb \Ga/\Ga_H}\;\;\sum_{(\ga_1,\ga_2)\in (\Ga_H\times \Ga_H)_\ga\bb(\Ga_H\times \Ga_H)}  f( h_1^{-1}\ga_1^{-1} \gamma \ga_2 h_2),\eeq
where the sum has only a finite number of nonzero terms.\\
Applying  (\ref{xigam}) and (\ref{sum}) to $|f|$ and using first the Fubini Theorem for positive functions and then for integrable functions, we obtain$$c_{\xi_\Ga,\xi_\Ga}(f)=\sum_{[\ga]\in \Ga_H\bb \Ga/\Ga_H} \int_{\Ga_H\bb H}\int_{\Ga_H\bb H}\sum_{(\ga_1,\ga_2)\in (\Ga_H\times \Ga_H)_\ga\bb (\Ga_H\times \Ga_H)}  f( h_1^{-1}\ga_1 \gamma \ga_2 h_2)dh_1\; dh_2$$
\ber\label{expr1}$$=\sum_{[\ga]\in \Ga_H\bb \Ga/\Ga_H} \int_{(\Ga_H\times\Ga_H)_\ga\bb H\times H}f( h_1^{-1}\gamma h_2)   d(h_1,h_2),$$\eer\\

the  integral and the sum being absolutely convergent.

As by assumption the group $(H\times H)_\ga$ is unimodular , using  the transitivity property of invariant measures on homogeneous spaces (see \cite{Be} Chap. II, \textsection 3), we have 
\ber\label{expr2}$$ \int_{(\Ga_H\times\Ga_H)_\ga\bb H\times H}f( h_1^{-1}\gamma h_2)  d(h_1,h_2)$$
$$=\int_{ (H\times H)_\ga\bb H\times H}\Big(\int_{(\Ga_H\times \Ga_H)_\ga\bb (H\times H)_\ga} f( h_1^{-1}u_1^{-1}\gamma u_2h_2) d(u_1,u_2)\Big) d(h_1,h_2)$$
$$ =vol((\Ga_H\times \Ga_H)_\ga\bb (H\times H)_\ga)\;\int_{ (H\times H)_\ga\bb H\times H}f( h_1^{-1}\gamma h_2)  dh_1\; dh_2.$$\eer
We deduce from this that the volume $vol((\Ga_H\times \Ga_H)_\ga\bb (H\times H)_\ga)$ is finite. Applying the above equality to $|f|$, we deduce that  the orbital  integral 
$I(\ga,f)$ of $f\in C_c(G)$ at $\ga$ is absolutely  converging. Thus we obtain the first assertion of the Proposition. \me

\no Therefore (\ref{expr1}) and (\ref{expr2}) give $$c_{\xi_\Ga,\xi_\Ga}(f)= \sum_{[\ga]\in \Ga_H\bb \Ga/\Ga_H}vol((\Ga_H\times \Ga_H)_\ga\bb (H\times H)_\ga)\;\int_{ (H\times H)_\ga\bb H\times H}f( h_1^{-1}\gamma h_2)  d(h_1,h_2).$$
Then the relative trace formula follows from (\ref{spect}).\qed
%%%%%%%%%%%%%%%%%%%%%%%%%%%%%%
%%%%%%%%%%%%%%%%%%%%%%%%%%%%
\section{Relative pseudocoefficient with  small support.}\label{PSC}
To define relative pseudocoefficient for relative discrete series, we need to review the abstract Plancherel formula. Here, we assume that $G$ is a real reductive group and $H$ is a spherical subgroup, in order to ensure finite multiplicities.\me

We denote by $\widehat{G}$ the unitary dual of $G$ and pick for every equivalence class $[\pi]$ a representative $(\pi,\cH_\pi)$. We keep notation of section \ref{SRTF}. \\
The abstract Plancherel formula Theorem for the spherical variety $\cZ:=G/H$ asserts the following: For every $[\pi]\in\hat{G}$, there exists a Hilbert space $\cM_{\overline{\pi}}\subset ( \overline{\cH}_\pi^{-\infty})^H$ (note that $\cM_\pi$ is finite dimensional, this induces a Hilbert space structure on $\textrm{Hom}(\cM_{\overline{\pi}},\cH_\pi)\simeq \cM_{\overline{\pi}}^*\otimes \cH_\pi$ and $\cM_\pi= \cM_{\overline{\pi}}^*\subset (\cH_\pi^{-\infty})^H$), such that the Fourier transform 
$$\cF:C_c^\infty(G)\to \int_{\hat{G}}^\oplus \textrm{Hom}(\cM_{\overline{\pi}},\cH_\pi) d\nu(\pi),$$
$$F\mapsto \cF(F)=\big(\cF(F)_\pi\big)_{\pi\in\hat{G}},\quad \cF(F)_\pi(\overline{\xi})=\overline{\pi}(F)\overline{\xi}\in\cH_\pi^\infty$$
extends to a unitary isomorphism. Here $\nu$ is a certain radon measure on $\hat{G}$ whose mesaure class is uniquely determined. The precise form of the measure depends on the chosen scalar products on the various stalks $\textrm{Hom}(\cM_{\overline{\pi}},\cH_\pi)$. We have 
$$\Vert F\Vert_{L^2(\cZ)}=\int_{\hat{G}} H_\pi(F) d\mu(\pi)$$
where $H_\pi$ are Hermitian forms which are defined as 
$$H_\pi(F)=\sum_{j=1}^{m_\pi}\Vert \overline{\pi}(F)\overline{\xi}_j\Vert^2_{\cH_\pi},\quad F\in C_c^\infty(\cZ),$$
for $\overline{\xi}_1,\ldots, \overline{\xi}_{m_\pi}$  an orthonomal basis of $\cM_{\overline{\pi}}$.
%%%%%%%%%%%%%%%%%%%%%%%%%%%%%%%%%
%%%%%%%%%%%%%%%%%%%%%%%%%%%%%%%%%%%%%%
\begin{defi}\label{pseudocoeff} Let $(\pi_0,\cH_{\pi_0})\in\hat{G}$ be a relative discrete series for $\cZ=G/H$, ie. which admits an embedding in $L^2(G/H)$. \begin{enumerate}\item  A function $f\in C^\infty_c(G/H)$ is a    relative pseudocoefficient for $\pi_0$ if
\begin{enumerate} \item for $\mu$-almost all $\pi\in \hat{G}$ distinct from $\pi_0$ and $\xi\in (\cH_\pi^{-\infty})^H$, then  $c_{\xi,\xi}(f)=0$,
\item there exists $\xi_0\in \cM_{\pi_0}$ such that $c_{\xi_0,\xi_0}(f)\neq 0$.
\end{enumerate}

\item Let $\xi_0\in (\cH_{\pi_0}^{-\infty})^H$.  A function $f\in C^\infty_c(G/H)$ is a  strong relative pseudocoefficient for $(\pi_0, \xi_0)$ if for any unitary irreducible representation $(\pi,\cH_\pi)$ of $G$ which admits an $H$-invariant distribution vector $\xi\in (\cH_\pi^\infty)^{',H}$, we have 
$$( \pi_{-\infty}(f) \xi,\xi)\neq 0\textrm{ if and only if } \pi=\pi_0\;\textrm{ and }\; \xi=c\xi_0 \textrm{ for some constant } c.$$
\end{enumerate}
\end{defi}
The relative trace formula for $\Ga\bb G$ allows to determine,  in some cases, if a relative discrete series $\pi$ for $G/H$ occurs in the spectral decomposition of $L^2(\Ga\bb G)$    and has a nonzero period (ie. $\xi_{\Ga,\pi}\neq 0$).
\begin{defi} We say that $f\in C^\infty_c(G/H)$ has small support relative to $\Ga$ if $I(f,\ga)\neq 0$ for $\ga\in \Ga$ implies that $\ga\in \Ga_H$.
\end{defi}
 Let us  assume that $G/H$ has a relative discrete series $(\pi_0,\cH_0)$. Then  $\cH_0$ can be realized as a subspace of $L^2(G/H)$ and the map $\xi_0 : \psi\in \cH_0^\infty \to \psi(1)$ is an $H$-invariant distribution vector. 
%%%%%%%%%%%%%%%%%%%%%%%%%%%%%
 \begin{prop}\label{period} 
  If  there exists a strong  relative pseudocoefficient $f$ for $(\pi_0, \xi_0)$ with  small support relative to $\Ga$ then $\pi_0$ occurs in $L^2(\Ga\bb G)$ with a nonzero period.
 \end{prop}
 \dem 
By definition, if $f$ is  a strong relative pseudocoefficient for $(\pi_0, \xi_0)$  with small support relative to $\Ga$, then the geometric side of the relative trace formula (\ref{RTF}) is reduced to the term $f(1)$ and the spectral side to the term $(\pi_0(f)\xi_{\Ga,\pi_0},\xi_{\Ga,\pi_0})$, hence we obtain the Proposition.\qed

We will precise the notion of  small support relative to $\Ga$ in the case of symmetric spaces. 
 We assume that $H$ is the fixed point group   of an involution $\si$ of  $G$. 
 \begin{lem}\label{Hcocompact} If $\Ga$ is a $\si$-stable cocompact discrete subgroup of $G$ then $\Ga_H=\Ga\cap H$ is a cocompact subgroup of $H$.
 \end{lem}
 \dem Let $(h_n)$ a sequence of $H$. As $\Ga\bb G$ is compact, extracting  possibly a subsequence of $(h_n)$, we can find a sequence $(\ga_n)$ in $\Ga$ such that $(\ga_n h_n)$ converges in $G$. Since $\Ga$ is $\si$-stable, the sequence of  $\ga_n\si(\ga_n)^{-1}=\ga_nh_n\si(\ga_nh_n)^{-1} $  is a converging sequence in $\Ga$. Hence, it is constant for $n$ large enough. Thus, there exists $n_0\in\N$ such that for $n\geq n_0$, we have $\ga_n\si(\ga_n)^{-1}=\ga\si(\ga)^{-1}$ where $\ga:=\ga_{n_0}\in\Ga$. This implies that $\ga^{-1}\ga_n\in\Ga\cap H$, 	and the sequence $(\ga^{-1}\ga_nh_n)$ converges. This proves that $\Ga\cap H\bb H$ is compact.\qed

 We assume that  $G$ is a reductive group in the Harish-Chandra class.  Let $\te$ be a Cartan involution of $G$ which commutes with $\si$. Then  $K:=G^\te$ is a maximal compact subgroup of $G$.\me

Let $\g g=\g k\oplus \g p=\g h\oplus \g q$ be the decomposition of the Lie algebra $\g g$ of $G$ in eigenspaces for $\te$ and $\si$ respectively.\me

We fix a maximal abelian subspace $\g a$ in $\g p\cap \g q$ and we denote by $A$ the analytic subgroup of $G$ with Lie algebra $\g a$. Then, we have the Cartan decompositions 
$$G=K\exp\g p=KAH.$$ 
We fix a $K$-invariant norm $\Vert\cdot \Vert$ on $\g p$ and we define a $K$-invariant function $\tau$ on $G$ by
$$\tau(k\exp X)=\Vert X\Vert,\quad k\in K, X\in \g p.$$

\no For $R>0$, let $A_R:=\{a\in A; \tau(a)<R\}$ be the ball of radius $R$ in $A$.

\no  We set  $$r_\Ga:=\inf_{g\in G, \ga\in \Ga-\{1\}}\tau(g^{-1} \ga g).$$ 
\begin{prop}\label{smallsup} Let $G$ and $\si$ be as above. Let $\Ga$ be a $\si$-stable cocompact discrete subgroup of $G$. Moreover, we assume that $\Ga$ is torsion-free. Then \begin{enumerate} \item  $r_\Ga>0$.
\item  Let $f\in C_c^\infty(G/H)$ be compactly supported  in   $KA_{r_\Ga/2}H$. Then $f$ has small support relative to $\Ga$.
\end{enumerate}\end{prop}
\dem  {\it 1.} This property is asserted in \cite{DGW}. We give a proof for sake of completeness.\\
If $r_\Ga=0$ then there would exist  two sequences  $(g_n)$ of $ G$ and $(\ga_n)$ of $ \Ga$
with $\ga_n\neq 1$ for all $n\in\N$, such that $\tau(g_n^{-1}\ga_ng_n)$ converges to $0$. Then 
\ber\label{lim}$g_n^{-1}\ga_ng_n=k_n \exp X_n$ with $k_n\in K$ and $X_n\in\g p$ with $\dis\lim_{n\to +\infty} \Vert X_n\Vert=0$.\eer\\
As $\Ga\bb G$ is compact, possibly changing $(g_n)$ and $(\ga_n)$ and extracting subsequences, we can assume that $(g_n)$ converges to $g\in G$ and $(k_n)$ converges to $k\in K$. Thus using (\ref{lim}), we see that the sequence $(\ga_n)$ converges, hence, as $\Ga$ is discrete,  it is constant and equal to $\ga:=\ga_{n_0}$ for $n\geq n_0$. Going to the limit in (\ref{lim}), we obtain $g^{-1}\ga g=k$. Therefore $\ga$ belongs to the discrete compact,  so finite,  group $gKg^{-1}\cap \Ga$. This implies that $\ga$ is a torsion element, thus $\ga=1$. This contradicts the hypothesis that all $\ga_n$ are distinct from $e$. This proves the first assertion.\me

\no{\it 2.}   Let $h_1, h_2$ in $H$ and $\ga\in\Ga$ such that $f(h_1^{-1}\ga h_2)\neq 0$. Then, the point $g:=h_{1}^{-1} \ga h_2$ belongs to $KA_{r_\Ga/2}H$. Therefore, we have 
 $g\si(g)^{-1}\in KA_{r_\Ga}K$, thus $\tau(g\si(g)^{-1})=\tau(h_1^{-1}\ga \si(\ga)^{-1} h_1)<r_\Ga$. By definition of $r_\Ga$, this implies that $\ga\si(\ga)^{-1}=e$, hence $\ga\in \Ga_H$. \qed
%%%%%%%%%%%%%%%%%%%%%%%%%
\section{Strong relative pseudocoefficient for some hyperbolic spaces.}\label{Strong}

The aim of this part is to construct strong relative pseudocoefficients associated to some  relative discrete series of  some hyperbolic spaces.

\subsection{Preliminaries.}
Let   $\K=\R,\C$ or $\H$ be the classical fiels of real, complex numbers or quaternions respectively.  Let $x\mapsto \bar{x}$ denotes the standard (anti-)involution of $\K$. Let  $p>2, q\geq 1$ be two integers. We consider the hermitian form $[\cdot,\cdot]$ on $\K^{p+q}$ given by
$$[x,y]=x_1\bar{y}_1+\ldots x_p\bar{y}_p-x_{p+1}\bar{y}_{p+1}-\dots -x_{p+1}\bar{y}_{p+q},\quad  (x,y)\in \K^{p+q}.$$

Let $G=U(p,q,\K)$ denote the group of $(p+q)\times (p+q)$ matrices preserving $[\cdot,\cdot]$. Let $H=U(1,\K)\times U(p-1,q, \K)$ be the stabilizer of $x_0= (1,0, \ldots, 0)\in \K^{p+q}$ in $G$. Then $H$ is the fixed point group of the involution   $\si$ of $G$  given by $\si(g)=JgJ$, where $J$ is the diagonal matrix with entries $(-1,1,\ldots,,1)$. The reductive symmetric space $G/H$ (of rank $1$) can be identified with the projective hyperbolic space $\X=\X(p,q,\K)$ (see \cite{AFJ} \textsection 2.) :
$$\X=\{z\in\K^{p+q}; [z,z]=1\}/\sim,$$\\
where $\sim$ is the equivalence relation $z\sim zu,\;\; u\in\K^*, |u|=1$.\me

The group $K=K_1\times K_2=U(p,\K)\times U(q,\K)\subset G$ is the maximal compact subgroup of $G$  consisting of elements fixed by  the classical Cartan involution   $\te$ of $G$, $\te(g)=(g^*)^{-1}$, which commutes to $\si$. Here $g^*$ denotes the conjugate transpose of $g$. \\
Recall that  $\g g=\g k\oplus \g p=\g h\oplus \g q$ are the decompositions of the Lie algebra $\g g$ of $G$ in eigenspaces for $\te$ and $\si$ respectively.\\
We define the one parameter abelian subgroup $A=\{a_t;\; t\in\R\}$ by 
$$a_t=\left(\begin{array}{ccc} \cosh t& 0 & \sinh t\\ 0 & I_{p+q-2}& 0\\ \sinh t & 0 & \cosh t\end{array}\right)$$
where $I_j$ denotes the identity matrix of size $j$. Then, the Lie algebra  $\g a$  of $A$ is a maximal abelian subspace of $\g p\cap \g q$.  Let  $W$ be the  Weyl group of $A$ in $G$. The nontrivial element of $W$ acts on $A$ by $a_t\mapsto a_{-t}$.\me

The Cartan decomposition $G=KAH$ holds and gives rise to the use of polar coordinates on X:
\beq\label{Cartan}(k,t)\in K\times \R^+\mapsto k a_tH,\eeq\\
and the map  $(k,t)\in K/K\cap M\times]0,+\infty[\mapsto k a_tH $ is a diffeomorphism on its image.\me

The centralizer $M$ of $A$ in $H$ is the subgroup of matrices 
$$\left(\begin{array}{ccc} u& 0& 0\\ 0&v&0\\ 0& 0&u\end{array}\right)$$
where $u\in\K^*$, $|u|=1$  and $v\in U(p-1,q-1,\K)$. \me

Hence, the homogeneous space $K/K\cap M$ can be identified with the projective image $\Y$ of the product of unit spheres $\S^p(\K)\times \S^q(\K)$:
\beq\label{K/KM}K/K\cap M\simeq \Y=\{ y\in\K^{p+q}; |y_1|^2+\ldots + |y_p|^2= |y_{p+1}|^2+\ldots + |y_{p+q}|^2=1\}/\sim.\eeq

Let $P$ be the subgroup of $G$ which stabilizes the $\K$-line generated by $\ga^0=(1,0,\ldots 0,1)$. Then $P$ is a maximal parabolic subgroup of $G$ whose  unipotent radical will be denoted by $N$, and we have $P=MAN$ (\cite{Fa} V.1).\me

Let $d:=\textrm{dim}_\R\K$. We  set $\dis \rho:=\frac{1}{2}(dq+dp)-1.$\me

We  recall some results about   spherical distributions of positive type on $\X$ given in \cite{Fa}. As   $\X$ is a symmetric space of rank  $1$,  the   algebra of left $G$-invariant differential operators on $\X$ is generated by the Laplace-Beltrami operator $\De$ corresponding to the natural pseudo-riemmannian structure. The  Laplace-Beltrami operator comes, up to a scalar, from the action of the Casimir of $\g g$ on $C_c^\infty(\X)$.\me

\no We denote by $\cD_{s,H}'(\X)$    the space of spherical distributions $\Te$ such that $\De\Te=(s^2-\rho^2)\Te.$\me
 
\no We keep notation of \textsection 1 for representations.

For $s\in\C$, we define the character $\de_s$ of $P$ by $\de_s(ma_tn)=e^{st}, m\in M, a_t\in A, n\in N$ and we denote by $(\pi_s:=ind_P^G \de_{s-\rho},\cH_s)$ the normalized induced representation.\\
According to (\cite{Fa}  page 395)
there exists a normalized $H$-invariant distribution vector $\xi_{s}\in( \cH_{-s}^{-\infty})^{H}$ (denoted by $u_s\in \Er_s'(\Xi)$ in loc. cit.), such that $s\to \xi_s(\phi) (\phi\in\cH_{-s}^\infty)$ is holomorphic on $\C$. We define the spherical distribution $U_s$ by
 (see \cite{Fa} Définition 5.1.)
 \beq\label{Us} U_s(f):=\langle (\pi_{-s})_{-\infty}(f) \xi_{s},\xi_{-s}\rangle,\quad f\in C_c^\infty(\X).\eeq\\
%%%%%%%%%%%%%%%%%%%%%%%%
 According to (\cite{Fa} Proposition 5.4  and Theorem 7.3), the spherical distribution $U_s$ satisfies the following properties:
 \ber\label{vpUs} \begin{enumerate}\item for $f\in C_c^\infty(\X)$, the map $s\to U_s(f)$ is holomorphic on $\C$,
 \item  $\De U_s=(s^2-\rho^2)U_s$,
\item $U_s=U_{-s}$ for $s\in \C$,
\end{enumerate}\eer

According to (\cite{Fa} IX, (2 a), Proposition 9.1, Théorème 9.2 and Proposition 9.3),   there exists at most a unique, up a scalar, eigendistribution   of positive type in $\cD'_{s,H}(\X)$  for $s\in\C$, except when $dq$ is even and $s=0$.

If $dq$ is odd, these distributions are the $U_s$ for $s\in i\R$, and the  $\ep_sU_s$  with $\ep_s=\pm1$ for a set of real $s$.

If $dq$ is even and $s\neq 0$, these distributions are the $U_s$ for $s$ in the union of $i\R$ and  some set of real $s$, and  the distributions denoted $(-1)^{r+1}\te_r\in \cD_{\rho+2r, H}'(\X)$ for $r\in \N$.

When $dq$ is even and $s=0$, denoting by $1$ the constant function equals to $1$, the distributions of positive type in $\cD'_{0,H}(\X)$ are given by $A(-\te_0)+B$ with $A,B\geq 0$.\me

Let $(\pi,V)$ be an irreducible unitary representation of $G$ and $\xi\in (V^{-\infty})^H$. Then the generalized matrix coefficient $c_{\xi,\xi}$ is a distribution of positive type. Moreover, as $\pi$ has infinitesimal character and as the Laplace-Beltrami operator comes from the Casimir of $\g g$, the distribution $c_{\xi,\xi}$ is an eigendistribution for the Laplace-Beltrami operator, hence, up a positive scalar, it is one of the  distributions  above.\me

Let us assume moreover that $(\pi,V)$ is a relative discrete series of $\X$, that is a subrepresentation of $L^2(\X)$. Let $\xi_V$ be the evaluation at $1$ of the elements of $V^\infty\subset C^\infty(\X)$. We say that the distribution $T:=c_{\xi_V,\xi_V}$ is associated to the  relative discrete series $(\pi,V)$ of $\X$.\me

By (\cite{Fa} Theorem 10 and page 432), up a positive scalar, the distributions associated to  relative discrete series of $\X$ are 
\begin{enumerate}\item if $\K=\R$ and $q$ is odd:\\
$$\ep_r U_{\rho+2r+1}\;\;\textrm{ for }\;\; r\in \Z \;\;\textrm{ such that }\;\; \rho+2r+1>0$$
where $\ep_r=(-1)^{r+1}$ if $r\geq 0$ and $\ep_r=1$ if $0<\rho+2r+1<2\rho$.
\item if $dq$ is even:\\
$$U_{\rho+2r}\;\;\textrm{ for }\;\; r\in -\N^* \;\;\textrm{ and }\;\; 0<\rho+2r<\rho,$$
\ber\label{RDS} and the $(-1)^r\te_r$ for $r\in\N$, which belong to $\cD'_{s_r,H}(\X)$ with $s_r=\rho+2r$.\eer
\end{enumerate}
\begin{prop}\label{nopseudo} Let  $(\pi,V)$ be a  relative discrete series of $\X$ whose associated distribution of positive type is of the form $\ep U_{s_0}$ for some $s_0\in\R$ and $\ep=\pm 1$. Then $(\pi,\xi_V)$ has no strong relative  pseudocoefficient.
\end{prop}
\dem Let $f\in C_c^\infty(\X)$ such that $\ep U_{s_0}(f)\neq 0$. By the holomorphy of $s\mapsto U_s(f)$ (see (\ref{vpUs})), the complex numbers  $U_s(f)$ for $s\in i\R$ are  not identically equal to $0$. This implies the Proposition.\qed

%%%%%%%%%%%%%%%%%%
%%%%%%%%%%%%%%%%%%%%%%%%%%%%%%%%%%%%%%%%%
%%%%%%%%%%%%%%%%%%%%%%%%%%%%%%%
\subsection{ $K$-types of  relative discrete series}\label{Ktype}
If this section, we assume that $\K=\C$ or $\H$. Hence, in particular $dq$ is even.  \me

We want to use the results of \cite{DFJ} where the groups are connected and semisimple.
For $\K=\C$, we have $G=Z G'$ where $Z:=\{z I_{p+q}; z\in\C^*, |z|=1\}$ is central and contained in $H\cap K$ and $G'$ is equal to $SU(p,q)$ which is semisimple and  connected. For $\K=\H$, then $G\simeq Sp(p,q)$ which is semisimple and connected (see \cite{Kn} Chap. I \textsection1).\me

For $r\in \N$, we will denote by $(\rho_r, V_r)$ the  relative discrete series whose associated distribution is $(-1)^r\te_r$, and by $\eta_r$ the element of $(V_r^{-\infty})^H$ such that $c_{\eta_r,\eta_r}=(-1)^r\te_r$. Note that $\eta_r$ is equal to a positive multiple of the evaluation at $1$ of the elements of $V_r^\infty\subset C^\infty(\X)$.\me

\no We will review the structure of $K$-module of $V_r$. For this, we introduce some notation.\me

\no By (\ref{K/KM}), the space $C^\infty(K/K\cap M)$ can be identified with the subspace of functions $f\in C^\infty (\S^p(\K)\times \S^q(\K))$ such that $f(u\zeta)=f(\zeta), \zeta\in \S^p(\K)\times \S^q(\K), u\in \K^*$ such that $|u|=1$.  According to (\cite{Fa} page 399), for $l,m\in\N$, we set 
$$\cY_{l,m}=\{f\in C^\infty(K/K\cap M); \De_1 f=-l(l+dp-2) f,\;\;  \De_2 f=-m(m+dq-2) f\}$$
where $\De_1$ and $\De_2$ are the Laplace-Beltrami operators in the spheres $\S^p(\K)$ and $\S^q(\K)$ respectively. 

Let $E$ be the set of elements $(l,m)\in\N\times \N$ such that $\cY_{l,m}\neq \{0\}$. By (\cite{Fa} page 399), $(l,0)\in E$ if and only if $l$ is even.\me

Let $r\in\N$. We set 
\beq\label{Er} E_r:=\{(l,m)\in E; l-m\geq dq+2r\}.\eeq\\
Then by (\cite{Fa} top of page 421),  the decomposition of $V_r$ as  $K$-module is given by
\beq\label{VrK}V_r=\bigoplus_{(l,m)\in E_r} \cY_{l,m}.\eeq\\

Let $l\in 2\N$ be even. We consider the function $\om_{l,0}$ defined in (\cite{Fa} bottom of page 406) for $\K=\C$ or $\H$. Let $\zeta=(\zeta_1,\ldots, \zeta_{p+q})$ be the coordinates on $\S^p(\K)\times \S^q(\K)$. Then $\om_{l,0}$ is given by:
\beq\label{oml0}  \om_{l,0}(\zeta)=\int_{U(1,\K)} F(\zeta_1u) du,\eeq\\
for some function $F$. \\
 Let us prove that 
\ber\label{ominv} $\om_{l,0}$ is biinvariant by $K\cap H$.\eer\\
By definition, $\om_{l,0}$ is right invariant by $K\cap M$. Recall that  $K=K_1\times K_2$ with $K_1=U(p,\K)$ and $K_2=U(q,\K)$. By loc. cit. top of page 407, we have  $\om_{l,0}\in\cY_{l,0}$, thus it is right invariant by $K_2$. But we have

$$K\cap M=\left\{\left(\begin{array}{cccc} u& 0& 0&0\\ 0&v_1&0&0\\ 0&0&v_2&0\\0&0& 0&u\end{array}\right), \quad u\in \K^*,|u|=1,  v_1\in U(p-1,\K), v_2\in U(q-1,\K)\right\}$$ and
$$K\cap H=\left\{\left(\begin{array}{ccc}  u& 0&0\\ 0&v_1&0\\0&0&w_2\end{array}\right), \quad  u\in \K^*,|u|=1,  v_1\in U(p-1,\K), w_2\in U(q,\K)\right\},$$
hence $K\cap H=(K\cap M)K_2$. Then $\om_{l,0}$ is right invariant by $K\cap H$.\me

 As $\K^*=\R^{+*}U(1,\K)$ where $\R^{+*}$ is central in the multiplicative group $\K^*$, we deduce from (\ref{oml0}) that $\om_{l,0}$  is left invariant by $K\cap H$. Hence (\ref{ominv}) follows. \me

\no We will determine the $K$-type of $\om_{l,0}$. 

Let $K'=K\cap SU(p,q)$ for $\K=\C$, or  $K'=K$ for $\K=\H$. We denote by $\g k'$  the Lie algebra of $K'$. Then $K$-types with $K\cap H$-fixed vectors coincides with $K'$-types with $K'\cap H$-fixed vectors since, for $\K=\C$, we have  $K=ZK'$ where $Z\subset K\cap H$ is central in $K$.

 We fix a maximal abelian subspace $\g t$ of $i(\g k'\cap \g q)$. As  $K'/K'\cap H$ is of rank $1$, the dimension of $\g t$ is  equal to $1$.  We choose a short positive root $\ga$ of  $\g t_\C$ in $\g k'_{\C}$. Then the roots of $\g t_\C$ in  $\g k'_{\C}$ are  of the form $\pm \ga,\pm 2\ga$. We identify $\C$ to $\g t_\C^*$ by the map $\la\mapsto \la\ga$. \\
 By the Cartan - Helgason Theorem, 
 \ber\label{phpoids}if $l\in\N$ is even, then $l$ is the highest weight a representation of $K$ with a unique, up to a scalar, nonzero $K\cap H$-fixed vector. \eer\\
  Let $\widehat{K} $  denotes the set of equivalence classes of unitary irreducible representations of $K$ and $(\widehat{K} )_{K\cap H}$ the subspace of that representations having a non-trivial $(K\cap H)$-fixed vector. For $\mu\in\widehat{K}$, let $\chi_\mu$  denotes its character and  $d_\mu$ its dimension. We set 
 $$\chi_\mu^H(k):=\int_{K\cap H} \chi_\mu(kh) dh,$$
where Haar measures on compact groups are normalized so that their volume are equal to $1$.  Then,
 \ber\label{chimuh} the function $\chi_{\mu}^H$ is, up to a scalar, the only function on $K$ of type $\mu$ which is biinvariant by $K\cap H$. \eer

 \begin{lem}\label{mur} Let $\mu\in (\widehat{K} )_{K\cap H}$ be the representation with highest weight $l\in2\N$. Then
 \begin{enumerate}\item $\mu $ is the unique $K$-type of $\cY_{l,0}$ having a nonzero  $K\cap H$-invariant vector.
 \item The multiplicity of $\mu$ in $\cY_{l,0}$ is equal to $1$ and $\om_{l,0}$ is contained in this $K$-type. 
  \item  $\om_{l,0}=C_l \chi_{\mu}^H$ with $C_l\neq 0$.
\end{enumerate}
 \end{lem}
 \dem Let $\mu'$ be a representation of $K$ contained in  $\cY_{l,0}$ and having a nonzero fixed vector by $K\cap H$. Then it has a highest weight of the form $k\ga$ where $k\in\N$ is even. The formula for the value of the Casimir operator acting on a highest weight representation implies that $k=l$, hence $\mu'=\mu$. \me

\no  As $\om_{l,0}$ is $K\cap H$-biinvariant by  (\ref{ominv}), we deduce easily {\it 2.} and  {\it 3.} from (\ref{chimuh}).\qed\\
We come back to the structure of $K$-module of the relative discrete series $(\rho_r,V_r), r\in\N$.\\
\ber\label{lrmur}For $r\in\N$, we denote by  $\mu_r\in (\widehat{K})_{K\cap H}$ the representation  with highest weight $l_r=dq+2r$ and we set 
$\om_{\mu_r}:=\om_{l_r,0}$. \eer\\
By (\ref{VrK}), we have $\om_{\mu_r}\in \cY_{l_r,0}\subset V_r$.

 \begin{lem}\label{Vr} Let $r\in\N$. There exists a unique  $K\cap H$- invariant function $\varphi_r$ in $V_r$ of   type $\mu_r$ such that $\varphi_r(1)=1$.\\
Moreover, there exists a constant $C'_r\neq 0$ such that $\varphi_r(k)=C'_r\chi_{\mu_r}^H(k)$ for $k\in K$.\end{lem}
\dem If $\varphi_r$ satisfies the first assertion of the Lemma then   the restriction of $\varphi_r$ to $K$ is a nonzero $K\cap H$- biinvariant function of type $\mu_r$.  Hence by (\ref{chimuh}), this restriction is proportional to $\chi_{\mu_r}^H$ and the second assertion follows.\me
 
Let us prove the first assertion. We first treat the case $\K=\H$. By (\cite{FJO} Table 2), there is at most one  relative discrete series for $L^2(\X)$ with a given eigenvalue of the Laplace-Beltrami operator. As $\mu_r$ satisfies (2.6) of loc. cit. (where $\langle \la,\la\rangle$ has to be replaced by  $\langle \al,\al\rangle$), Theorem 2.2 in loc. cit. implies that $V_r$ contains a unique, up to a scalar, $K\cap H$ invariant function  of type $\mu_r$ denoted there by $\psi_\la$ with $\la=\rho+2r$. By definition, it satisfies $\psi_\la(1)=\psi^0_\la(1)$ where $\psi^0_\la$ is given in loc. cit. (2.5). The formula defining $\psi^0_\la$ shows that $\psi^0_\la(1)\neq 0$. Thus the function $\varphi_r:=\psi_\la/\psi_\la(1)$ satisfies the first assertion of the proposition.\me

For $\K=\C$, we proceed similarly by first going through the quotient by the center $Z\subset K\cap H$.\qed\me

\no Let $\mu\in (\widehat{K} )_{K\cap H}$ with highest weight $l\in2\N$.  For  $s\in\C$,  the vector  $(\pi_{-s})_{-\infty}(\chi_{\mu })\xi_s=(\pi_{-s})_{-\infty}(\chi_{\mu} ^H)\xi_s$ is an analytic vector for $\pi_s$. Thus, using Lemma \ref{mur} {\it 3.}, we can define 
\beq\label{Usom} \ga_{l }(s):=\langle (\pi_{-s})_{-\infty}(\om_{l,0})\xi_s, \xi_{-s}\rangle= U_s(\om_{l,0}).\eeq

An explicit expression of  $U_s(\Phi)$ for a $K$-finite function $\Phi$ in $C_c^\infty(\X)$ is obtained in  (\cite{Fa} page 407). This expression   allows us to calculate the function $\ga_{l}(s)$ in the next Lemma, which  is given for granted in \cite{DFJ}.

%%%%%%%%%%%%%%%%

\begin{lem}\label{Usom} Let $\mu\in (\widehat{K} )_{K\cap H}$ be the representation with highest weight $l\in2\N$. Let $s\in\C$. According to (\cite{Fa} page 405), we define the function  $$\be_{l ,0}(s)=b_{l }\frac{(s-\rho)(s-\rho-2)\ldots (s-\rho-l_r+2)}{\Ga((s-\rho+l +dp)/2)},\quad\textrm { where } b_{l } \textrm{ is a nonzero constant}.$$
Then, we have
$$ \ga_{l }(s)=U_s(\om_{l,0})= c'_{l} \be_{l,0}(s)\be_{l,0}(-s)$$
for some  constant $c'_l>0$.
\end{lem}
\proof We consider the function $A(t)=   (e^t+e^{-t})^{dp-1}(e^t-e^{-t})^{dq-1}$ according to (\cite{Fa} page 403).
 We can find  a sequence $(F_n)_{n\in\N}$ of $C_c^\infty(]0,+\infty[)$ such that $\textrm{supp}(F_n)\subset ]\dfrac{1}{2n}, \dfrac{1}{n}[$ and 
$$\int_0^{+\infty} F_n(t) A(t) dt=1.$$\\
Therefore, we have\beq\label{dirac} \lim_{n\to +\infty}\int_0^{+\infty} F_n(t) A(t) g(t) dt=g(0),\quad \textrm{ for } \; g\in C_c^\infty(\R).\eeq\\
Using the Cartan decomposition (\ref{Cartan}),  we define the function $\Phi_n$ on $\X$  by $\Phi_n(ka_tx^0)=\om_{l,0}(k) F_n(t),\; k\in K, t\in [0,+\infty[.$ Hence, each $\Phi_n$ is of  type $\mu$ and belongs to $C_c^\infty(\X)$. By  (\cite{Fa} page 407), we have 
$$U_s(\Phi_n)=c_{l}   \be_{l,0}(s)\be_{l,0}(-s)\int_0^{+\infty} \Psi_{l,0}(t,s) F_n(t) A(t) dt \times \int_{K/K\cap M} \om_{l,0}(k)^2 dk,$$
where $c_{l}$ is a nonzero constant and $\Psi_{l,0}$ is given in term of the hypergeometric function by
$$ \Psi_{l,0}(t,s) =(\cosh t)^{s-\rho}\; _2F_1\big(\frac{\rho-s+l}{2}, \frac{\rho-s-dp+2-l}{2}, \frac{dq}{2}, \tanh ^2t).$$\\
%%%%%%%%%%%%%%%%%%%%%%%%%%%%%
Since $\Psi_{l,0}(0,s)=1$, we deduce from   (\ref{dirac})   that 
$$U_s(\om_{l,0})=c'_{l}   \be_{l,0}(s)\be_{l,0}(-s)
$$
with $c'_{l}=c_{l}  \int_{K/K\cap M} \om_{l,0}(k)^2 dk$. Thus, we obtain the Lemma.\qed

 %%%%%%%%%%%%%%%%%%%%%%%%%%%%%%

%%%%%%%%%%%%%%%%%%%%%%%%

 %%%%%%%%%%%%%%%%%%%%%%%%%%%%%%%%%%
\subsection{Existence of strong relative pseudocoefficients for certain relative discrete series.}
In this section, we assume that $dq$ is even and $\K=\C$ or $\H$.\me

Existence of strong relative pseudocoefficient for relative discrete series $(\rho_r,V_r), r\in \N$ is an easy consequence of the next Lemma.  This  Lemma  corresponds to  (\cite{DFJ} Lemma 9), but the proof given in loc. cit. is slightly incomplete. We will give here a more precise and modified proof.\me

Recall that $A=\{a_t; t\in \R\}$ is the abelian subgroup of $G$ corresponding to the  maximal abelian subspace $\g a$ of $\g p\cap \g q$. We fix $R>0$ and we denote by $A_R:=\{a_t; |t|<R\}$  the ball  of radius $R$ in $A$.
\begin{lem}\label{lem9} Let $r\in\N$ and $\mu_r$ be the $K$-type of highest weight $l_r=dq+2r$. Then, there exists a function $f\in C_c^\infty(\X)$  with  support in $KA_RH$, sum of a $K$-invariant function and of a function of type $\mu_r$, such that 
\begin{enumerate}

\item $\te_{r }(f)= 1$,
\item $\te_{r'}(f)=0$ for  $r'\in \N$ such that  $r'>r$,
\item $U_s(f)=0$ for $ s\in\C$.
\end{enumerate}
\end{lem}

\dem  We first recall some results of \cite{DFJ} on the Paley-Wiener space of $\X$. Notice that for $\K=\C$, we have to go through the quotient by the center $Z\subset K\cap H$ to apply these results.
\me

Let  $(\pi,V)$ a unitary irreducible representation of $G$ in a Hilbert space $V$. If $\mu\in \widehat{K}$ then its contragredient is equivalent to $\mu$ and $P_\mu:=\pi_{-\infty}(d_\mu \chi_\mu)$  is well defined as the projection of $V^{-\infty}$ onto the $\mu$-isotypic component $V_\mu\subset V^\infty$. If $\xi$ is an $H$-invariant distribution vector, then 
\beq\label{proj} P_\mu\xi=d_\mu\pi_{-\infty}(\chi_{\mu}^H)\xi\in V_\mu^{K\cap H},\eeq
hence $P_\mu\xi=0$ if $\mu\notin (\widehat{K} )_{K\cap H}$.\me

\no Let $\g a_\C^*$ be the complexification of the dual $\g a^*$ of $\g a$. Recall that $W$ denotes the Weyl group of $A$ in $G$, hence we can consider the action of $W$ on $\g a$.\me

Let $R>0$. Let $PW(\g a)_R$ denotes the  space of entire functions $\Psi$ on $\g a^*_\C$ which   satisfy
$$ \forall N\in\N,\;\; \sup_{\la\in \g a^*_\C} (1+\Vert\la\Vert)^N e^{-R \Vert \textrm{Im}\;\la\Vert}|\Psi(\la)|<+\infty.$$
Then the classical Paley-Wiener space $PW(\g a)$ is the union of $PW_R(\g a)$ for $R\in]0,+\infty[$.

Restricting to the $W$-invariant functions, the classical Fourier transform is a bijection of $C_c^\infty(\g a)^W$ to $PW(\g a)^W.$ \me

According to (\cite{DFJ} Theorem 1 and its remark), we have the following result:
\ber\label{PW} Let   $\Psi\in PW_R(\g a)^W$. Let $\mu\in \widehat{K}_{K\cap H}$. Then, there exists a unique function $f\in C_c^\infty(G/H)$ of  type $\mu$, supported  in $KA_RH$ such that the following holds:

For all  unitary irreducible representation $(\pi,V)$ of $G$  and for all $H$-invariant distribution vector $\xi\in V^{-\infty}$  such that $\pi_{-\infty}(\De)\xi= (\la^2-\rho^2)\xi$, we have 
$$\pi_{-\infty}(f) \xi=\Psi(i\la) P_\mu\xi.$$
\eer\\
%%%%%%%%%%%%%%%%%%%%
We recall that  $r\in \N$ is  fixed and $l_r=dq+2r$ is the highest weight of $\mu_r$ (see (\ref{lrmur})).  We set  $s_r:=\rho+2r$. As   $c_{\eta_r,\eta_r}=(-1)^{r+1}\te_r\in \cD'_{s_r,H}(\X)$ (see (\ref{RDS})), we have  $\rho_r(\De)\eta_r=(s_r^2-\rho^2)\eta_r$. Let us prove the following result:
\ber\label{teG} Let  $G\in PW_R(\g a)^W$ such that $G(is_r)\neq 0$. 
Then, there exists  $f_1\in C_c^\infty(G/H)$  of  type $\mu_r$  supported in $KA_RH$   such that \begin{enumerate}
\item $\te_{r }(f_1)= 1$,
\item $\te_{r'}(f_1)=0$ for  $r'\in \N$ such that  $r'>r$,
\item $U_s(f_1)=C\ga_{l_r}(s) G(is),\; s\in\C$,
for some nonzero constant $C$.\end{enumerate}
\eer
  We   apply (\ref{PW}) to $\mu=\mu_r$. Then there exists $g_1\in C^\infty_c(\X)$ of type $\mu_r$ supported in $KA_RH$ such that
\beq\label{fr}(\rho_{r})_{-\infty}(g_1)\eta_r=G(is_r) P_{\mu_r} \eta_r\eeq

and\beq\label{fs}(\pi_{-s})_{-\infty}(g_1)\xi_s=G(is) P_{\mu_r}\xi_s,\;\textrm{for } s\in \C.\eeq

By Lemma \ref{mur} and (\ref{proj}), we have $P_{\mu_r} \eta_r=d_{\mu_r}(\rho_r)_{-\infty}(\chi_{\mu_r}^H)\eta_r$ and the analytic vector $(\rho_r)_{-\infty}(\chi_{\mu_r}^H)\eta_r$ is a $K\cap H$-invariant function of type $\mu_r$ in $V_r$. By Lemma \ref{Vr}, we obtain that  $(\rho_r)_{-\infty}(\chi_{\mu_r}^H)\eta_r=C' \varphi_r$ for some constant $C'$. \\Let us prove that $C'\neq 0$. By Lemma \ref{Vr} again, the function  $\varphi_r$ coincides with $C'_r\chi_{\mu_r}^H$ on $K$. Since $\rho_r(\chi_{\mu_r}^H)$ is, up to a scalar,  the projection  on the $K\cap H$-fixed vectors in $V_r$, we have $\rho_{r }(\chi_{\mu_{r }}^H) \varphi_{r }=\varphi_r$. Thus we deduce $$C'\langle \varphi_r,\varphi_r\rangle=\langle(\rho_{r })_{-\infty}(\chi_{\mu_{r }}^H)\eta_{r },\varphi_{r }\rangle=\langle \eta_{r }, \rho_{r }(\chi_{\mu_{r }}^H) \varphi_{r }\rangle=\Big(\rho_{r }(\chi_{\mu_{r }}^H) \varphi_{r }\Big)(1)=\varphi_r(1)=1,$$
hence $C'\neq 0$.  Then $\te_r(g_1)=C' G(is_r)\neq 0$. We set 
$$f_1:=\frac{g_1}{C' G(is_r)}.$$
Then  $f_1$ satisfies $$\te_r(f_1)=1.$$
and  we obtain the first assertion of (\ref{teG}). \me 

Let  $r'\in\N$ such that $r'>r$. As $f_1$ is of type $\mu_r$ and $\mu_r$ is not a $K$-type of    $V_{r'}$ by (\ref{VrK}), we have $\te_{r'}(f_1)=0$.\me

To prove the last assertion of (\ref{teG}), we consider the property (\ref{fs}) for  $s\in\C$.  By (\ref{proj}), we have  $P_{\mu_r}\xi_s=d_{\mu_r}(\pi_{-s})_{-\infty}(\chi_{\mu_r})\xi_s=d_{\mu_r}(\pi_{-s})_{-\infty}(\chi_{\mu_r}^H)\xi_s$. Using Lemma \ref{mur} and Lemma \ref{Usom}, this leads to 
$$U_s(f_1)=\frac{d_{\mu_r}}{C_r C' G(is_r)}G(is) U_s(\om_{\mu_r})=\frac{d_{\mu_r}}{C_r C' G(is_r)}G(is) \ga_{l_r}(s).$$\\ 
This finishes the proof of (\ref{teG}).\me

 By definition (see Lemma \ref{Usom}), we have  $$\ga_{l_r}(s)= c'_{l_r} \be_{l_r, 0}(s)\be_{l_r,0}(-s),$$
and $$\be_{l_r, 0}(s)= b_{l_r}\frac{(s-\rho)(s-\rho-2)\ldots (s-\rho-l_r+2)}{\Ga((s-\rho+l_r+dp)/2)},$$ \\
where $l_r=dq+2r$.\\
Recall that $\rho=\dfrac{1}{2}(dp+dq)-1$, then we have $-\rho+l_r+dp= \rho+2r+2= s_r+2$ with $s_r=\rho+2r$. Thus,  we can write $$\Ga((s-\rho+l_r+dp)/2)=  2^{-l_r}(s+s_r)(s+s_r-2)\ldots (s-\rho+dp)\Ga((s-\rho+dp)/2).$$
Hence we obtain 
$$ \ga_{l_r}(s)=\frac{P(s)}{(s+s_r)(-s+s_r)Q(s)}\ga_0(s)$$
where $$P(s)= 2^{l_r}c'_{l_r}b_{l_r}^2(s-\rho)(s-\rho-2)\ldots (s-\rho-(l_r-2))(-s-\rho)(-s-\rho-2)\ldots (-s-\rho-(l_r-2))$$
 and $$Q(s)=(s+s_r-2)(s+s_r-4)\ldots (s-\rho+dp)(-s+s_r-2)(-s+s_r-4)\ldots (-s-\rho+dp).$$
By assumption we have $dq\geq 2$, hence $2r\leq dq+2r-2=l_r-2$. Then, $P_1(s):=\dfrac{P(s)}{(s+s_r)(-s+s_r)}$ and $Q(s)$   are even polynomials such that    $Q(s_r)\neq0$ and $P_1(s_r)\neq 0$. \me

By ((\ref{teG}) {\it 3.}), we obtain

$$U_s(f_1)=C G(is)\frac{P_1(s)}{Q(s)}\ga_0(s),\quad \textrm{with } \ga_0(s)=c'_0\frac{1}{\Ga(\frac{s-\rho+dp}{2})\Ga(\frac{-s-\rho+dp}{2})}.$$\\
Since $Q$ is even and $Q(  s_r)\neq 0$, we may choose an invariant differential operator $D\in\D(\X)$ such that 
$$U_s(Df_1)=\frac{Q(s)}{Q(s_r)} U_s(f_1),\quad s\in\C,$$\\
and $$\te_{r'}(Df_1)=\frac{Q(s_{r'})}{Q(s_r)} \te_{r'}(f_1),\quad r'\in\N.$$\\
Therefore the function $f_2:=Df_1$ satisfies the assertions {\it 1.} and {\it 2.} of the Lemma and 
$$U_s(f_2)=C \frac{G(is) P_1(s)}{Q(s_r)}\ga_0(s).$$

Since $PW_R(\g a)^W$ is stable by multiplication by an even polynomial function, 
applying  (\ref{PW}) to the  trivial $K$-type, we can find a $K$-invariant function $f_3\in C_c^{\infty}(\X)$ supported in $KA_RH$ such that $$U_s(f_3)=C\frac{G(is) P(s)}{Q(s_r)}\ga_0(s)=U_s(f_2).$$  
By (\ref{VrK}), the trivial representation is not a $K$-type of $V_{r'}$ for $r'\geq 0$, thus we have $\te_{r'}(f_3)=0$ for $r'\geq 0$. Therefore, the function $f=f_2-f_3$ satisfies the properties of the Lemma.\qed
%%%%%%%%%%%%%%%%%%
%THEOREM
\begin{theo}\label{Result} Let $R>0$ and $A_R:=\{a_t\in A; |t|<R\}$ be the ball of radius $R$ in $A$. Then, for all $r\in\N$, there  exists   a strong relative pseudocoefficient $f\in C_c^\infty(\X)$  supported  in $KA_RH$ for $(\rho_r,\eta_r)$. This means that  the function $f$   satisfies
$$\te_r(f)=1,\quad \te_{r'}(f)=0 \textrm{ for } r'\in\N,\;\; r'\neq r,$$
and 
$$ U_s(f)=0, \textrm{ for } s\in\C$$ \end{theo}
\dem  Lemma \ref{lem9} gives the result for $r=0$. Let  $r>0 $ and assume we have a strong relative pseudocoefficient $f_{r'}$ for $(\rho_{r'},\eta_{r'})$ ($r'< r$) supported in $KA_RH$. We denote by $\Psi_{r}$  the function  obtained  in  Lemma \ref{lem9}.  Then, the function  $\dis f_r=\Psi_r-\sum_{r'=0}^{r-1} \te_{r'}(\Psi_{r}) f_{r'}$ is a strong relative pseudocoefficient for $(\rho_{r},\eta_{r})$ with support in $KA_RH$.\qed

\section{Existence of $\si$-stable torsion free  cocompact  discrete subgroups of $U(p,q,\K)$.}\label{exis}
This section is entirely due to the kind help of R. Beuzart-Plessis and J. P. Labesse.\me

Let $\F$ be  a totally real number field of degree $[\F:\Q]=r>1$ and let $V_\infty=\{v:\F\hookrightarrow \R\}$ denote the finite set of real places of $\F$.

We consider the group $\bG'$ defined in  (\cite{SW} page 372), which depends on  $\K=\R,\C$ or $\H$. By loc. cit., there exists a unique archimedean place $v_1\in V_\infty$ such that
\beq\label{grpeR} \bG'(\F_{v_1})=\left\{\begin{array}{ll} SO(p,q)& \textrm{ for } \K=\R\\
SU(p,q)& \textrm{ for } \K=\C\\
Sp(p,q)& \textrm{ for } \K=\H \end{array}\right. \eeq\\
and   for $v\in V_\infty, \; v\neq v_1$,
\beq\label{gpev}\bG'(\F_v)=\left\{\begin{array}{ll} SO(p+q)& \textrm{ for } \K=\R\\
SU(p+q)& \textrm{ for } \K=\C\\
Sp(p+q)& \textrm{ for } \K=\H \end{array}\right.\eeq
Let $\bG=\textrm{Res}_{\F/\Q}\bG'$ be the group obtained by restriction of scalars. For  $v\in V_\infty, v\neq v_1$, we have  $\bG(\Q)=\bG'(\F)\subset \bG'(\F_v)$ and $\bG'(\F_v)$ is compact, hence each element of $\bG(\Q)$ is semisimple. Thus  the group $\bG(\Q)$ is anisotropic. \me
  
 Let $\A_\F$ and $\A_\Q$ be the ring of adeles of $\F$ and $\Q$ respectively. We denote by $\A_{\F,f}$ and $\A_{\Q,f}$ the subrings of finite adeles in $\A_\F$ and $\A_\Q$ respectively. 

Then, $\bG(\Q)$ is diagonally embedded in  $\bG(\A_{\Q,f})$ and by (\cite{PR} Theorem 5.5 (1)), the quotient $\bG(\Q)\bb \bG(\A_\Q)$ is compact. \me

\no Let $\tau$ be the  rational involution  of $\bG'$ denoted $\tau_{1,0}$ in (\cite{SW} \textsection 2.2). The involution of $\bG'(\F_{v_1})$ induced by $\tau$, again denoted by $\tau$,  is simply the restriction to $\bG'(\F_{v_1})$  of the involution $\si$ of $U(p,q,\K)$ defined in section 3. 
The involution $\tau$ defines a continuous automorphism of $\bG(\A_\Q)=\bG'(\A_\F)$ preserving $\bG(\A_{\Q,f})=\bG'(\A_{\F,f})$.\me
 
 Let $K_f$ be an open compact subgroup of $\bG(\A_{\Q, f})$. Since  $\bG(\Q)$ is discrete in $\bG(\A_\Q)$, 
the subgroup $\Ga(K_f):=\bG(\R)\cap \big( \bG(\Q) K_f )$ is a discrete subgroup of $\bG(\R)$. As  $\bG(\Q)$ is diagonally embedded in  $\bG(\A_{\Q,f})$, we have  also $\Ga(K_f)=\bG(\Q) \cap K_f$.
 
We have \ber\label{Gacomp} $\Ga(K_f)\bb \bG(\R)$ is compact.\eer

Let us give a proof for sake of completeness. We consider the map $\psi$ from $\bG(\R)K_f$ to $\Ga(K_f)\bb \bG(\R)$ given by $\psi(gk)=\Ga(K_f)g$. If $gk= \ga g' k'$ with $g,g'\in \bG(\R)$, $k,k'\in K_f$ and $\ga \in \bG(\Q)$, we have $g= \ga g' k' k^{-1}$. As $\bG(\A_{\Q,f})$ and $\bG(\R)$ commute, we obtain $g=\ga  k' k^{-1}g'$, hence $gg'^{-1}= \ga  k' k^{-1}\in \Ga(K_f)$. Therefore $\psi$ goes through the quotient $\bG(\Q)\bb \bG(\R)K_f/K_f$ and the induced map is a surjection from $\bG(\Q)\bb \bG(\R)K_f/K_f$ to $\Ga(K_f)\bb \bG(\R)$. \me

Thus it remains to prove that $\bG(\Q)\bb \bG(\R)K_f/K_f$ is compact.\me

The group $\bG(\R)K_f$   is an open subgroup of $\bG(\A_\Q)$. As $\bG(\Q)\bb \bG(\A_\Q)$ is compact, the number of $(\bG(\Q), K_f)$ double cosets, which are open, is finite. Then each of them is also closed, hence  compact. We deduce that $\bG(\Q)\bb \bG(\R)K_f/K_f$ is compact and  (\ref{Gacomp}) follows.%%%%%%%%%%%%
\begin{lem} Let $K_f$ be a compact open subgroup of $\bG(\A_{\Q,f})$. Then there exists a compact open subgroup $K'_f\subset K_f$ of $\bG(\A_{\Q,f})$ such that $\Ga(K'_f)$ is a torsion free subgroup of $\Ga(K_f)$ of finite index. 
\end{lem}
\dem We fix an embedding of $\bG$ in $GL(n)$ defined over $\Q$. We have 
$$GL(n,\Z)=GL(n,\Q)\bigcap\prod_{p\;\textrm{prime}}GL(n,\Z_p).$$

By the proof of (\cite{Bor} Proposition 2.2), the group $GL(n,\Z)$ contains a torsion free subgroup $\tilde{\Ga}$ of finite index of the form $\tilde{\Ga}=GL(n,\Q)\cap \tilde{K}_f$, where $\tilde{K}_f$ is a compact open subgroup of $\dis \prod_{p\;\textrm{prime}}GL(n,\Z_p)$.

Let $K'_f:=\tilde{K}_f\cap K_f\subset \bG(\A_{\Q,f})$. Then $\Ga(K'_f)=\bG(\Q)\cap K'_f\subset \tilde{\Ga}$ is without torsion.  As $K'_f$ is a compact open subgroup of $K_f$, it is of finite index. It follows that $\Ga(K'_f)$ is a finite index in $\Ga(K_f)$. \qed
\begin{cor}\label{gama} For each open compact subgroup $K_f$ of $\bG(\A_{\Q,f})$, there exists a $\tau$-stable open compact subgroup $K^0_f\subset K_f$ such that  $\Ga(K^0_f)$ is a $\tau$-stable  torsion free cocompact discrete subgroup of $\bG(\R)$. 
\end{cor}
\dem  Let  $K'_f$ be the subgroup obtained in the previous Lemma. As $\tau$ is a continuous involution of $\bG(\A_\Q)$ which preserves $\bG(\A_{\Q,f})$, the subgroup $K^0_f:=\tau(K'_f)\cap K'_f$ is a $\tau$-stable compact open subgroup of $K_f$.  It is clear that  $\Ga(K^0_f)$ is $\tau$-stable since $K^0_f$ is $\tau$-stable. The properties of $K'_f$ implies that  $\Ga(K^0_f)$ is a torsion free cocompact discrete subgroup of $\bG(\R)$.  Hence we obtain the Corollary.\qed

%%%%%%%%%%%%%%%%%%%%%%
\begin{lem}\label{prod} Let $G=G_1\times G_2$ be the product of two locally compact groups with $G_2$  compact. Let $\Ga$ be a torsion free cocompact discrete  subgroup of $G$. Then the projection $\Ga_1$ of  $\Ga$ to $G_1$ is a torsion free cocompact discrete subgroup of $G_1$.
\end{lem}
\dem  If  $\Ga_1$ was not discrete, there would exists a sequence $(\ga_{1,n})$  of distinct elements of  $\Ga_1$ converging to a limit $l$. There exists  a sequence $(\ga_{2,n})$ in $G_2$ such that $\ga_n=(\ga_{1,n},\ga_{2,n})$ belongs to $\Ga$. Since $G_2$ is compact, extracting a subsequence, we can assume that $(\ga_{2,n})$ converges. Then the sequence $(\ga_n)$ convergences, hence it is constant for $n$ large enough as $\Ga$ is discrete. This implies that $(\ga_{1,n})$ is constant for $n$ large enough, which  contradicts the fact that the $\ga_{1,n}$ are distinct. Thus $\Ga_1$ is discrete.\me

Let us show that $\Ga_1$ is a cocompact subgroup of $G_1$. Let $(g_n)$ be a sequence in $G_1$. Since $\Ga$ is cocompact in $G$,  there exist a subsequence $(g'_n)$ of  $(g_n)$ and a sequence $(\ga_n)$ in $\Ga$ such that $(\ga_ng'_n)$ converges. Writing $\ga_n=(\ga_{1,n}, \ga_{2,n})$ with $\ga_{i,n}\in G_i$, we deduce  that $(\ga_{1,n} g'_n)$ converges, hence $\Ga_1$ is cocompact in $G_1$.\me

Let $\ga_1\in \Ga_1$ and $r\in \N^*$ such that $\ga_1^r=1$. Let $\ga_2\in G_2$ such that $(\ga_1,\ga_2)\in \Ga$. Then the sequence $\big((\ga_1,\ga_2)^n\big)=(\ga_1^n,\ga_2^n)$ remains in a compact set, thus it admits  a converging subsequence $((\ga_1,\ga_2)^{k(n)})$. As $\Ga$ is discrete, this subsequence is constant for $n$ large enough. Hence $(\ga_2^{k(n)})$ is constant for $n$ large enough. This implies that there exists $s\in\N^*$ such that $\ga_2^s=1$. For $m$ a multiple of $r$ and $s$, we have $(\ga_1,\ga_2)^m=1$. As $\Ga$ is torsion free, this leads to $\ga_1=\ga_2=1$, hence $\Ga_1$ is torsion free.\qed

\begin{lem}\label{compact} If $\Ga$ is cocompact in $SO(p,q)$ (resp., in $SU(p,q)$) then $\Ga$ is cocompact in $O(p,q)$ (resp. in $U(p,q)$).\end{lem}
\dem 
This follows from the fact that $SO(p,q)$ (resp., $SU(p,q)$) is cocompact in $O(p,q)$ (resp., $U(p,q)$).\qed\me

\no By (\ref{grpeR}) and (\ref{gpev}), there is a compact group $\Om_{ \K}$, depending on   $\K$, such that 
\beq\label{gr}\bG(\R)=\left\{\begin{array}{ll} SO(p,q)\times \Om_{\R} & \textrm{ for } \K=\R\\
SU(p,q)\times \Om_{\C} & \textrm{ for } \K=\C\\
Sp(p,q)\times \Om_{\H} & \textrm{ for } \K=\H\end{array}\right. .\eeq
We denote by $G_{1,\K}$ the first factor of this decomposition.

\begin{prop}\label{taustable}  For sufficiently small $\tau$-stable  open compact subgroup $K_f$ of $\bG(\A_{\Q, f})$, the projection $\Ga_1(K_f)$ of  $\Ga(K_f)=\bG(\Q)\cap K_f$  onto  $G_{1,\K}$ according the decomposition (\ref{gr}) is a $\si$-stable torsion free cocompact discrete subgroup of $U(p,q,\K)$.
\end{prop}
\dem By Corollary \ref{gama}, we can choose $K_f$  sufficiently small so that $\Ga(K_f)$ is a $\tau$-stable torsion free cocompact discrete subgroup of $\bG(\R)$. 
 By Lemma \ref{prod}, the subgroup $\Ga_1(K_f)$ is a $\tau$-stable torsion free cocompact discrete subgroup of $G_{1,\K}$. Since the involution $\tau$ coincides with $\si$ on $G_{1,\K}$,  the Proposition follows from    Lemma \ref{compact}. \qed

\begin{theo} Let $\K=\C$ or $\H$.  Let $K_f$  and $\Ga_1(K_f)$ be as in  Proposition \ref{taustable}. Then the relative discrete serie $(\rho_r,V_r)$ of $U(p,q,\K)$ occurs with a nonzero period in $L^2(\Ga_1(K_f)\bb U(p,q,\K))$.\end{theo}
\dem To apply the relative trace formula (\ref{RTF}), we have to verify that $\Ga_1(K_f)$ and $H$ satisfy assumptions (\ref{hyp}). The group $H=U(1,\K)\times U(p-1,q,\K)$ is unimodular. By (\ref{grpeR}) and (\ref{gpev}), each element of $\bG(\Q)$  is semisimple, hence each element $\ga$ of $\Ga_1(K_f)$ is semisimple. By (\cite{V} Part II, chap. 2 Proposition 13) , the centralizer of $\ga\si(\ga)^{-1}$ in $G=U(p,q,\K)$ is reductive since $\ga\si(\ga)^{-1}$ is semisimple. Moerover, this centralizer is $\si$-stable, hence the centralizer $Z_H(\ga\si(\ga)^{-1})$ of $\ga\si(\ga)^{-1}$ in $H$ is reductive. As the identity map induces an isomorphism  from  $(H\times H)_\ga$ to $Z_H(\ga\si(\ga)^{-1})\times H$, we deduce that $(H\times H)_\ga$ is reductive, hence unimodular. The quotient $(\Ga_1(K_F)\cap H)\bb H$ is compact by Lemma \ref{Hcocompact}. Therefore the assumptions (\ref{hyp}) are all satisfied. By Theorem \ref{Result} there exists a strong relative pseudocoefficient for $(\rho_r,V_r)$, with arbitrary small support.  As $\Ga_1(K_f)$ is torsion free, Proposition \ref{smallsup} and Proposition \ref{period} give the result.\qed.

\section{ Non existence of  relative pseudocoefficients for $G(\C)/G(\R)$.}\label{GC}
Let  $G$ be a   connected, simply connected  complex semisimple Lie group.  Let  $H$ be a real form of   $G$ and $\si$ be the conjugation of $G$ relative to $H$. We denote by $\g g$ and $\g h$  the Lie algebras of $G $ and $H$ respectively. Let   $\g g=\g h\oplus \g q$ be the decomposition of $\g g$ relative to $\si$. Hence we have $\g q=i\g h$. Recall that a Cartan subspace of $\g q$ is a maximal abelian subspace made of semsimple elements. Then the map $\g a\to i\g a$ is an isomorphism from the set of Cartan subalgebras of $\g h$ to the set of Cartan subspaces of $\g q$ which preserves $H$-conjugacy classes.\me

By (\cite{OM2} Theorems 1 and 2), the symmetric space $G/H$ has relative discrete series if and only if $\g q$ has a compact Cartan subspace, or equivalently, if $\g h$ is a split real form of $\g g$. The goal of this section is to establish that there is no  relative pseudocoefficient for any relative discrete series of $G/H$ (see Theorem \ref{nopseudo} below). This result will follow from the inversion formula  of orbital integrals (see \cite{Ha3} Théorème 6.15). \me

We assume that $H$ is split and we fix a split Cartan subalgebra $\g a_0$ of $\g h$. Let $\Ga_{\g a_0}$ be the lattice of $X\in a_0$ satisfying $\exp 2iX=1$ and let  $\Ga_{\g a_0}^*$ be its dual lattice.  Let $P=LN$ be a $\si$-stable Borel subgroup of $G$ with Levi subgroup $L= \exp(\g a_0+i\g a_0)$. For $\mu\in\Ga_{\g a_0}^*$, we define the character $\de_\mu$ of $P$ by $\de_\mu(\exp(X+iY)n)=e^{i\mu(Y)}$ for $X,Y\in\g a_0$ and $n\in N$. We denote by $(\pi_\mu, \cH_\mu)$ the normalized induced representation $(ind_P^G \chi_\mu, \cH_\mu)$.

 Then by (eg.  \cite{De} Proposition 5 and \cite{Ha1} \textsection 3. Application 1.), the relative discrete series of $G/H$ is given by $(\pi_\mu, \cH_\mu)$ where  $ \mu\in \Ga_{\g a_0}^*$ is regular and the space of $H$-invariant distribution vectors of $\pi_\mu$ is generated by $\xi_\mu$ defined by the integration over $H/H\cap P$ 
$$\xi_\mu(\psi)=\int_{H/H\cap P} \psi(h) d\dot{h}, \quad \psi\in \cH_\mu^\infty,$$
where $d\dot{h}$ is a semi-invariant measure on $H/H\cap P$.\me

We denote by $\Ga_{\g a_0, reg}^*$ the set of regular elements in $\Ga_{\g a_0}^*$ and by $\Ga_{\g a_0, sing}^*$ its complementary in $\Ga_{\g a_0}^*$.\me

\no We   recall some facts about regular elements in $G/H$ and orbital integrals  (see \cite{OM1} and \cite{Ha3} \textsection 1 and \textsection 2).

Let  $\varphi$ be the map from $G/H$ to $G$ defined by $\varphi(gH)=g\si(g)^{-1}$.
A  semisimple element $x\in G/H$ is  regular if $\varphi(x)$ is semisimple and regular in $G$ in the usual sense. Let $(G/H)_{reg}$ be the open dense subset of semisimple regular elements of $G/H$. To  a Cartan subalgebra  $\g a$  of $\g h$, we associate the Cartan subspace $A$ of $G/H$ consisting of elements $x$ such that $\varphi(x)$ centralize $\g a$. If $x\in (G/H)_{reg}$  then the centralizer $\g a:=Z_{\g h}(\varphi(x))$ of $\varphi(x)$ in $\g h$ is a Cartan subalgebra of $\g h$ and $x$ belongs to the Cartan subspace $A$ associated to $\g a$. If $Y\in \g h$ then $(\exp iY)H$ is regular in $G/H$ if and only if $Y$ is regular in $\g h$ in the usual sense.\\
 If $V$ is a subset or a subalgebra of $\g h$ or of its dual space $\g h^*$, then $V_{reg}$ will denote the set of regular elements in $V$.\me

Let $D_G$ be the Weyl discriminant of $G$ as a complex semisimple  Lie group. The orbital integral $\cM(f)$ of $f\in C_c^\infty(G)$ is the function $\cM(f)\in C^\infty((G/H)_{reg})$ defined by 
$$\cM(f)(x)=|D_G(\varphi(x))|^{1/2}\int_{H/Z_H(\g a)}f(h\cdot x) dh,$$
where $\g a=Z_{\g h}(\varphi(x))$ and $Z_H(\g a)$ is the centralizer of $\g a$ in $H$.
\no By (\cite{Bou} \textsection 8), 
\ber\label{mfsup} for all Cartan subset $A$ of $G/H$, there exists a compact subset $U\subset A$ such that for all $x\in A_{reg}-U$, then $\cM(f)(x)=0$.\eer\me

By (\cite{Bou} \textsection 8), orbital integrals are characterized by three properties and the above condition on its support. According to the definition of  (\cite{Ha3} \textsection 2), orbital functions are   $H$-invariant functions in $C^\infty((G/H)_{reg})$ which sastify the same properties as orbital integrals except the condition on the support (\ref{mfsup}).

%%%%%%%%%%%%%%%%%%%%%%%%%%
%%%%%%%%%%%%%%%%%%%%%%%%%%
\begin{theo}\label{nopseudo}Let $\mu_0\in \Ga_{\g a_0^*}$ be regular. Let $(\pi_{\mu_0},\cH_{\mu_0})$ be the relative discrete serie associated to $\mu_0$. Then there exists no relative pseudocoefficient for $(\pi_{\mu_0},\cH_{\mu_0})$.
\end{theo}
\dem 
We recall the inversion formula of orbital integrals (see \cite{Ha3} Théorème 6.15).\\
Let $\te$ be a Cartan involution of $\g h$ commuting with $\si$. We fix a system $[Car(\g h)]$ of $\te$-stable  representants  of the $H$-conjugacy classes of Cartan subalgebras of $\g h$. We may and will assume that $\g a_0\in [Car(\g h)]$.

Let $\g a\in [Car \g h]$.  Let $\g a=\g a_I\oplus \g a_R$ be its decomposition with respect to $\te$.  We denote by $\Ga_{\g a}$ the  lattice made of elements $X\in \g a_R$ such that $\exp 2iX=1$ and by $\Ga_{\g a}^*$ its dual lattice.
For $M=G$ or $H$, we denote by $W_M(\g a)$  the quotient of the normalizer  of $\g a$ in $M$ by the centralizer   of $\g a$ in $M$ and we set $\cW_{\g a}:=W_G(\g a)/W_H(\g a)$. We fix a positive system $\psi$ of imaginary roots of $\g a_\C$ in $\g g$. Let $\g a_{I,reg}^*$ be the set of $\la\in\g a_I^*$ such that $\langle \la,\al\rangle\neq 0$ for all $\al\in \psi$. Then, for $\g a=\g a_0$, the set $\cW_{\g a_0}$ is reduced to the trivial element which we denote by $1$ and $\psi=\emptyset$.

By (\cite{Ha2} Théorème 6.1) and (\cite{Ha3} Théorème 5.3), to each $(\la, y)\in \Ga_{\g a}^*+\g a_{I,reg}^*\times\cW_{\g a}$, we can associate an $H$-invariant eigendistribution $\Te(\la,y,\psi)$ and an orbital function $F(\la,y,\psi)$ such that, for all $f\in C_c^\infty(G/H)$, we have 

\beq\label{InvF}\cM(f)=\sum_{\g a\in [Car \g h]} c_{\g a} \sum_{\mu\in \Ga_{\g a}^*} \int_{\g a_I^*} \sum_{w\in \cW_{\g a}^{-1}} \sum_{y\in\cW_{\g a}}F(w(\mu+\la), y,\psi) \langle \Te(-w(\mu+\la), y,\psi),f\rangle d\la,\eeq
where $c_{\g a} $ for $\g a\in[Car \g h]$  are constants depending only on the choices of measures.\me

Moreover, for  $\g a=\g a_0$ and $\mu\in \Ga_{\g a_0, reg}^*$ then   $\Te(-\mu, 1,\emptyset)$ is equal to the generalized matrix coefficient $c_{\xi_{\mu},\xi_\mu}$ associated to $(\pi_\mu,\xi_\mu)$. \me

Let $\mu_0\in\Ga^*_{\g a_0,reg}$. Assume that there exists a   relative pseudocoefficient $f\in C_c^\infty(G/H)$ for $(\pi_{\mu_0},\cH_{\mu_0})$. Let us prove that 
\ber\label{quest}the function $F_f$ defined by the right hand side of  (\ref{InvF}) do not satisfy (\ref{mfsup}).\eer  This contradiction will achieve the proof of the Theorem.  For proving this, we will evaluate  $F_f$ at $(\exp iX)H$ for some $X\in \g h$ regular.\me

Let $\g a\in [Car \g h]$, $y,w\in\cW_{\g a}$ and $\la\in \Ga^*_{\g a}+\g a^*_{I,reg}$. By (\cite{Ha3} Théorème 5.8), we have $F(w\la, y,\psi)((\exp iX)H)\neq 0$ if and only if $y=1$. By (\cite{Ha2} Proposition 6.4), there exists a unitary irreducible representation $(\pi_\la, \cH_\la)$ and an $H$-invariant distribution vector $\xi^w_\la\in \cH_\la^{-\infty}$ such that the generalized matrix coefficient $c_{\xi^w_\la,\xi^w_\la}$ is equal to $(i)^{|\psi|} \ep_{w,\la,\psi}\Te(-w^{-1}\la, 1,\psi)$ where $\ep_{\la,w,\psi}=\pm 1$. Since  $(\pi_\la, \cH_\la)$  belongs to the support of the Plancherel formula of $G/H$ (see \cite{Ha3} Théorème 7.4) and    $f$ is a  relative pseudocoefficient for $(\pi_{\mu_0},\xi_{\mu_0})$,  if $\g a\neq \g a_0$ or if $\g a=\g a_0$ and if $\la\in \Ga_{\g a_0,reg}^*$ is  different from $\mu_0$,  then we have  $\langle\Te(-w^{-1}\la, 1,\psi),f\rangle=0$.
Thus we obtain 
\beq\label{egal}F_f((\exp iX)H)=c_{\g a_0} \big(F(\mu_0)((\exp iX)H)+\sum_{\mu\in \Ga_{\g a_0,sing}^*} F(\mu) ((\exp iX)H) \langle \Te_\mu, f\rangle),\eeq
where $F(\mu):=F(\mu,1, \emptyset)$.\me

Since the Cartan subset $\exp(i\g a_0)H$ associated to $\g a_0$ is compact, we will evaluate $F_f$ on another Cartan subset. For this, we recall briefly the construction of $F(\mu)$   for convenience of the reader (see \cite{Ha3} \textsection 4).\\
The functions $F(\mu)$ for $\mu\in\Ga_{\g a_0}^*$ are obtained from the Fourier transform $\hat{\be}_{H\cdot X}$ of the Liouville measure $\be_{H\cdot X}$ on the orbit $H.X$ for $X\in \g h_{reg}$. \me

Let  $\g a$ be Cartan subalgebra of $\g h$. We fix $x\in G$ such that $x\cdot\g a_{0,\C}=\g a_\C$. By the properties of the Fourier transform of orbits (see \cite{V} Théorème I.7.7), there exist  constants $c(w,\cC,\cF)$  depending on connected components $\cF$ and $\cC$ of  $\g a_{_0,reg}^*$ and  $\g a_{reg}$ respectively, and on $w\in W_G(\g a_{0,\C})$,  such that,  we have 
$$\hat{\beta}_{H\cdot X}(\mu)\mid \textrm{det}(\textrm{ad}\;\mu)_{\g h^*/\g a_0^*}\mid^{1/2}=\sum_{w\in W_G(\g a_{0 \C})} c(w,\cC,\cF)e^{i\langle xw\mu, X\rangle}, \quad X\in\cC,  \mu\in\cF.$$
Moreover, 
\ber\label{cste}$c(w,\cC,\cF)\neq 0$ if and only if  $\textrm{Im}\langle xw\mu,X\rangle\geq 0$ for all  $X\in \cC$ and $\mu\in\cF$.\eer\\
 For $\ep>0$ and $U$ an open neighborhood of $0$ in the center of $\g h$, we denote by $\cV_{\ep,U}$    the  open set of $X\in U\oplus[\g h,\g h]$ such that the real part $\textrm{ Re}(\la)$ of each eigenvalue $\la$ of $\textrm{ad} X$ satisfies $|\textrm{ Re}(\la)|<\ep$. 
Then for $U$ and $\ep$ small enough, the set $\cW_{\ep,U}:=\bigcup_{h\in H} h\exp i\cV_{\ep,U}H$ is an open neighborhood of $1$ in $G/H$.\me

Let  $\mu\in \Ga_{\g a_{0}}^*$. We consider the set $\cF(\mu)$ of connected components $\cF$ of $\g a_{0,reg}^*$ such that $\mu$ belongs to the closure of $\cF$. Then for any connected component $\cC$ of $\g a_{reg}$, we have (see \cite{Ha3} \textsection 4)
\beq\label{fmu} F(\mu)((\exp iX) H)=\frac{1}{|\cF(\mu)|}\sum_{\cF\in\cF(\mu)}\Big(\sum_{w\in W_G(\g a_{0 \C})} c(w,\cC,\cF)e^{i\langle xw\mu, X\rangle}\Big),\quad X\in \cC\cap \cV_{\ep,U}.\eeq\\
Notice that if $\mu$ is regular, we have 
$$ F(\mu)(\exp iX H)=\hat{\beta}_{H\cdot X}(\mu)\mid \textrm{det}(\textrm{ad}\;\mu)_{\g h^*/\g a_0^*}\mid^{1/2},\quad X\in ( \cV_{\ep,U})_{reg}.$$
Let  $\al$ be a (real) root of  $\g a_0$ in $\g h$. Let $X_{\pm\al}\in\g h$  be root vectors in $\g h$ such that $H_\al:=[X_{-\al},X_{\al}]$ is the coroot of $\al$. Then $\g a_{\al}:=\R (X_\al-X_{-\al})\oplus \textrm{Ker}\; \al$ is a Cartan subalgebra of $\g h$ and $c_\al.\g a_{0\;\C}=\g a_{\al\;\C}$ where $c_\al:=\textrm{Ad}(\exp\;-i\frac{\pi}{4}(X_\al+X_{-\al}))$ is the usual Cayley transform. Then the imaginary roots of $\g a_\al$ in $\g h$ are $\be=c_\al(\al)$ and $-\beta$ and we have  $\g a_{\al}=\R iH_\beta + \textrm{ker}\;\beta$ with $ \textrm{Ker}\;\beta= \textrm{Ker}\;\al$. 
\me

We may and will choose   $Y\in \textrm{ Ker }\be\cap \cV_{\ep,U}$  such that for all   $t>0$, the element  $X_t:=itH_{\beta} +Y$ is regular in $\g a_\al$.  Let $\cC$ be the connected component  of $\g a_{\al, reg}$ which contains $X_t$, $t>0$. \\
Since the roots of  $\g a_\al$ take imaginary values on  $i\R H_{\beta}$ and real values on $ \textrm{ Ker }\be$ and $Y\in \textrm{ Ker }\be\cap \cV_{\ep,U}$, we have  $X_t\in \g a_\al\cap \cV_{\ep,U}$   for all $t>0$. For  $w\in W_G(\g a_{0,\C})$, we have  $ \langle c_\al w\mu, X_t\rangle= it\langle w\mu,H_\al\rangle+\langle w\mu,Y\rangle$ since $c_\al(Y)=Y$ and $c_\al H_\al=H_\be$. For $\mu\in \Ga_{\g a_0}^*$,   we obtain by (\ref{fmu})
$$F(\mu)((\exp iX_t)H)=\frac{1}{|\cF(\mu)|}\sum_{\cF\in\cF(\mu)}\Big(\sum_{w\in  W_G(\g a_{0\Bbb C})}c(w,\cC,\cF)e^{-t\langle w\mu,H_\al\rangle} e^{i\langle w\mu,Y\rangle}\Big),\quad\textrm{for all } t>0.$$\\
Using  (\ref{egal}), we deduce that $F_f(\exp iX_tH)$ is nonzero for all $t>0$. 
Therefore, since the set of $(\exp iX_t)H$ for $t>0$ is not included in a compact subset of the Cartan subset associated to $\g a_\al$, we deduce  the assertion (\ref{quest}).\qed

%%%%%%%%%%%%%%%%%%%%%%%%%%%%%%

\noindent P.~Delorme, Aix-Marseille Université, CNRS, Centrale Marseille, I2M, UMR 7373, 13453
Marseille, France.\\
{\it E-mail address}:  patrick.delorme@univ-amu.fr\me

\noindent P.~Harinck, CMLS, \'Ecole polytechnique, CNRS-UMR 7640, Université Paris-Saclay, Route de Saclay,  91128 Palaiseau Cedex, France.\\  {\it E-mail address}: pascale.harinck@ polytechnique.edu\me

\end{document}